\documentclass[11pt,reqno,a4paper]{amsart}
\usepackage{euscript,verbatim}
\newcommand{\pp}{\hspace{1.2em}}
\def\CC{\mathbb{C}}
\def\RR{\mathbb{R}}

\def\MM{\mathbb{M}}
\def\NN{\mathbb{N}}

\def\e{\epsilon}
\def\II{\mathbb{I}}
\def\I{\mathbb{I}}

\def\etal{\textit{et al.}}

\newtheorem{theorem}{Theorem}[section]
\newtheorem{lemma}[theorem]{Lemma}
\newtheorem{proposition}[theorem]{Proposition}
\newtheorem{corollary}[theorem]{Corollary}

\theoremstyle{definition} 

\newtheorem{example}[theorem]{Example}
\newtheorem{definition}[theorem]{Definition}

\theoremstyle{remark}
\newtheorem{remark}[theorem]{Remark}

\date{22 october 2002} 

\title[Products of Matrices and Multiperiodic Functions]{Products of non-stationary random matrices and Multiperiodic equations of several scaling factors}

\author[A.H. Fan]{ Ai-Hua FAN} 
\address{Ai-Hua FAN: \ CNRS UMR 6140 -- LAMFA, Universit\'e de Picardie, 80039 Amiens, France}
\email{ai-hua.fan@u-picardie.fr}
\urladdr{http://www.mathinfo.u-picardie.fr/fan/}

\author[B. Saussol]{Beno\^{\i}t SAUSSOL}
\address{Beno\^{\i}t SAUSSOL: \ CNRS UMR 6140 -- LAMFA, Universit\'e de Picardie, 80039 Amiens, France}
\email{benoit.saussol@u-picardie.fr}
\urladdr{http://www.mathinfo.u-picardie.fr/saussol/}

\author[J. Schmeling]{J\"org SCHMELING}
\address{J\"{o}rg SCHMELING:
     \ CNRS UMR 6140, LAMFA, Universit\'e de Picardie, 80039 Amiens, France  \&
  Department of Mathematics, University of Lund, P.O. Box 118, SE-221 00 LUND, Sweden}
\email{joerg@maths.lth.se}

\thanks{\textbf{Acknowledgement.} This work was partially done during the first author's visit to Lund University, Sweden
and the last author's visit to WIPM of the academy of sciences of China and 
 Wuhan University, China.}

\begin{document}
\begin{abstract}

 Let $\beta>1$ be a real number
 and $M: \mathbb{R}\rightarrow {\rm GL(\CC^d)}$ be a uniformly almost
periodic matrix-valued  function.
We study the asymptotic behavior of the product
$$
P_n(x) =M(\beta^{n-1}x)\cdots M(\beta x) M(x).
$$
 Under some condition we prove a theorem of Furstenberg-Kesten type
for such products of non-stationary random matrices. Theorems of
Kingman and Oseledec type are also proved.
The obtained results are applied to multiplicative functions defined
by commensurable scaling factors.
We get a positive answer to a Strichartz conjecture on the asymptotic
behavior of such multiperiodic functions.  The case where $\beta$ is
a Pisot--Vijayaraghavan number is well studied.

\end{abstract}

%\subjclass{Primary:28A80,42A38.}

\maketitle

%%%%%%%%%%%%%%%%%%

%%%%%%%%%%%%%%%%%%

\vspace{2em}

\section{Introduction}

Kingman's subadditive ergodic theorem was originally proved in
1968 \cite{Ki1,Ki2}. A more recent proof was given by Katznelson
and Weiss in 1982 \cite{KW}. It is one of the most important
results in ergodic theory. In this paper we consider the
following set-up which resembles a dynamical system without
invariant measure and try to get results similar  to Kingman's
theorem. Let $\beta
>1$ be a positive real number. Let $\{f_n\}$ be a sequence of
uniformly almost periodic functions (i.e. in the sense of Bohr, see Section~\ref{sec:2.1})
defined on the real line $\mathbb{R}$. Suppose the following
subadditivity condition is fulfilled
\[
     f_{n+m}(x) \le f_n(x) + f_m(\beta^n x)\ \
     \mbox{for a.e.} \ x \ \mbox{and all} \ n, m.
\]
where $a.e.$ refers to the Lebesgue measure. We would like to
study the almost everywhere convergence of $n^{-1} f_n(x)$. The
Kingman theorem applies in the special case where $\beta>1$ is an
integer and the $f_n$'s are periodic. The typical case in our mind is
\begin{equation}\label{log-product-matrix}
f_n(x) = \log \|M(\beta^{n-1} x) \cdots M(\beta x) M(x)\|
\end{equation}
where $M: \mathbb{R} \rightarrow \mbox{GL}(\mathbb{C}^d)$ is a
matrix-valued uniformly almost periodic function. We will prove
that the limit $\lim_{n \to \infty} n^{-1}f_n(x)$ exists almost
everywhere (a.e. for short) with respect to the Lebesgue measure
under the condition that the $n^{-1} f_n(x)$ have \emph{joint periods} (see
Theorem~\ref{king}). As a consequence, an Oseledec type theorem is proved
for the matrix products involved in \eqref{log-product-matrix}
(see Theorem~\ref{ose}). It is proved that the condition on the existence
of joint periods is satisfied when $\beta$ is a PV-number
(see Section~\ref{sec:3}).

 \pp

Our consideration is partially motivated by the study of
multiperiodic functions, already investigated by Strichartz \etal\
\cite{JRS}, Fan and Lau \cite{FL}, and Fan \cite{F}.
 By  a {\em Multiperiodic function}  of one real variable we mean any
function $F: \mathbb{R} \rightarrow \mathbb{R}$ which is a solution of a
functional equation of the  following form $$
      F(\xi) = f_1\left( \frac{\xi}{\rho_1}\right) F\left( \frac{\xi}{\rho_1}\right)
+ \cdots + f_d\left( \frac{\xi}{\rho_d}\right) F\left(
\frac{\xi}{\rho_d}\right)
      $$ where $d\geq 1$ is an integer;
      $\rho_1>1, \ldots, \rho_d>1$ are $d$ real numbers, called {\em scaling
      factors}; $f_1, \ldots, f_d$ are $d$ complex valued functions defined
      on the real line, called  {\em determining functions}. The equation will be called
      a {\em multiperiodic equation}.

\pp
      We will assume that the determining functions $f_j$ are
      periodic or almost periodic in the sense of Bohr, as is the case in most
      applications. We will also assume that the scaling factors
      $\rho_j$ are commensurable in the sense that $\rho_j$ are powers of some
      real number $\beta >1$. Without lost of generality, we assume that
      $\rho_j = \beta^j$ for $1\leq j \leq d$. Then the multiperiodic equation
      becomes
       \begin{equation}~\label{eq-multiperiodic}
      F(\xi) = f_1\left( \frac{\xi}{\beta}\right) F\left( \frac{\xi}{\beta}\right)
+ \cdots + f_d\left( \frac{\xi}{\beta^d}\right) F\left(
\frac{\xi}{\beta^d}\right)
      \end{equation}
      As far as we know, there is  few work done for the non--commensurable case
      which is much more difficult.

\pp In the literature, the case where $d=1$ and $\beta = 2$ (or an
arbitrary integer) has been studied, especially in the theory of
wavelets \cite{D}. In fact, the scaling function $\varphi$ of a
wavelet satisfies a scaling equation $$
     \varphi(x) = \sum a_n \varphi(2 x -n).
$$ The Fourier transform of $\varphi$ satisfies a multiperiodic equation
of the form \eqref{eq-multiperiodic} with only one scaling factor
$\beta =2$ and only one determining function $f_1(x)=f(x) =
\frac{1}{2} \sum_n a_n e^{inx}$.

\pp The scaling functions in wavelets constitute a class of
functions sharing a kind of similarity. More generally,
multiperiodic functions arise as  Fourier transforms of
self-similar objects such as Bernoulli convolution measures
($d=1$, $\beta>1$ being a real number and $f$ being a
trigonometric polynomial), inhomogeneous Cantor measures ($d$ may
be greater than $1$) or more general self-similar measures produced
by iterated function systems (see \cite{S}).

\pp In the case of one scaling factor (i.e. $d=1$, then we write
$f_1(x) = f(x)$), the existence of the solution of the multiperiodic
equation \eqref{eq-multiperiodic} is simple and is assured by
the consistency condition $f(0) =1$ and a regularity condition,
say $f$ is Lipschitz continuous. Actually the solution can be written as an
infinite product $$
     F(x) = \prod_{n=1}^\infty f \left(\frac{x}{\beta^n}\right).
$$ \pp
For the existence of the general equation
\eqref{eq-multiperiodic}, we have

\medskip

\noindent {\bf Theorem A.} \ {\it Let $d\geq 1$. Suppose  that the
determining functions $f_1, \cdots, f_d$ are Lipschitz continuous, and satisfy
the consistency condition $$
        f_1(0) + \cdots + f_d(0) = 1.
$$
Suppose furthermore  that $f_1(0) , \ldots, f_d(0) \in [0,+\infty)$.
Then equation~\eqref{eq-multiperiodic} admits a unique continuous
solution $F$ such that $F(0) =1$. }
\medskip

The proof of this theorem is postponed to Section~\ref{sec4.2}.

\pp
 Our study of equation \eqref{eq-multiperiodic} is converted to that of vector valued
equations of the form
\begin{equation}~\label{eq-vector}
     G(x) = M\left( \frac{x}{\beta}\right)  G\left( \frac{x}{\beta} \right)
\end{equation}
where  $M:  \mathbb{R} \rightarrow \mathcal{M}_{d \times d}(\mathbb{C}) $
is a matrix valued determining function and  $G : \mathbb{R} \rightarrow
\mathbb{R}^d$ is a  vector valued
 unknown function.  Matrix products will be
 involved in the study of equation \eqref{eq-vector}, which produces some difficulties.
 However,  equation~\eqref{eq-vector} is a simple recursive relation  because it contains only one
 scaling factor. Equation~\eqref{eq-multiperiodic} is equivalent to equation \eqref{eq-vector} with $M(x)$ and $G(x)$ equal respectively to
 \begin{equation}\label{Matrix-multip}
 \left( \begin{array}{ccccc}
       f_1(x) & f_2(\frac{x}{\beta}) & \cdots &  f_{d-1}(\frac{x}{\beta^{d-2}})& f_d(\frac{x}{\beta^{d-1}})\\
       1      & 0            & \cdots & 0                 & 0 \\
       0       &     \ddots   &  \ddots&\vdots                  & \vdots\\
      \vdots        &  \ddots             &         &  0                     & 0\\
     0         &   \cdots            &   0      &  1                     & 0
                            \end{array}
\right)
\text{and}\left( \begin{array}{c}
F(x)\\ F(\frac{x}{\beta}) \\ \vdots \\ F(\frac{x}{\beta^{d-1}})
\end{array} \right)
 \end{equation}

\pp We would like to know the asymptotic behavior at infinity
of the solution $G$. This is a natural question because $G$ often
represents Fourier transform of a function (a measure or a
distribution) and the asymptotic behavior at infinity
describes quantitatively the regularity of the solution.
Unfortunately, there is no closed form formula for $G$ in general
and the behavior of $G$ is rather complicated, as is shown by the
Fourier transform of the Cantor measure ($d=1$,  $\beta=3$ and $f(\xi)
= \cos  \xi$).

\pp Following \cite{JRS}, we will study the pointwise asymptotic
behaviors of $$
                 h_n (x) :=\frac{1}{n} \log |F(\beta^n x)|
$$ as $n \rightarrow \infty$.

We will prove that, under some conditions, the limit  $ \lim_{n
\rightarrow \infty}h_n(x)$ exists and is equal to a constant
almost everywhere with respect to Lebesgue measure.

This answers partially a questions in \cite{JRS}. More precisely,
we have the following results, whose proofs are postponed to
Section~\ref{sec4.3}.

\medskip

\noindent {\bf Theorem B.} \ {\it Suppose that the conditions in
Theorem A are satisfied. Furthermore,  suppose that the
determining functions $f_1, \cdots, f_d$ are either identically zero or
strictly positive and
$1$-periodic, and that $\beta>1$ is a Pisot number.
Let $F$ be the
solution of the equation~\eqref{eq-multiperiodic}.
 Then there is a constant $\mathcal{L}$ such that $$
      \lim_{n \rightarrow \infty} \frac{1}{n} F(\beta^{n} x) = \mathcal{L}
      \quad a.e.
$$}

 The constant $\mathcal{L}$ in the theorem is the leading Liapunov exponent
of the matrix $M(x)$ above, defined by
$$
                \mathcal{L}(M) = \inf_{n\geq 1} \frac{1}{n} \mathbb{M} \log \|
                M(\beta^{n-1}x) M(\beta^{n-2}x) \cdots M(\beta x) M(x)\|.
 $$
where $\mathbb{M} f$ denotes the Bohr mean of an almost periodic function
$f$ (see Section~\ref{sec:2.1} below).

\medskip

\noindent {\bf Theorem C.} \ {\it Suppose  that the conditions in Theorem A
are satisfied. Furthermore,  suppose that the determining functions $f_1,
\cdots, f_d$ are $1$-periodic, Lipschitz, and that $\beta>1$ is a Pisot
number with maximal conjugate of modulus $\rho$. Let $F$ be the solution of
the equation~\eqref{eq-multiperiodic}.
If
\[
\sup_x\frac{(1+|f_1(x)|+\cdots+|f_{d-1}(x/\beta^{d-2})|)
(|f_1(x)|+\cdots+|f_d(x/\beta^{d-1})|)}{|f_d(x/\beta^{d-1})|}
<\rho^{-1}
\]
then there is a constant $\lambda\in\RR$ such that
\[
\lim_{n\to\infty} \frac{1}{n} \log \sum_{j=0}^{d-1}|F(\beta^{n-j} x)| = \lambda
\quad\text{a.e.}
\]
}

In this case we do not know if the constant $\lambda$ equals the
leading Lyapunov
exponent of the matrix.

Theorem A will be proved in Section 4.2 as a special case of
 a more general result
(Theorem \ref{existence-G}), Theorm B and Theorem C in Section 4.3.
Both  Theorm B and Theorem C are consequences of our Kingman's Theorem (Theorem \ref{king})
and Oseledec's Theorem (Theorem \ref{ose}) which are discussed in Section 2.
In Section 3, we prove that the joint period condition
required in both Kingman's Theorem and Oseledec Theorem is satisfied
when $\beta>1$ is a Pisot number.

\bigskip

\section{Kingman theorem and Oseledec theorem}
\setcounter{equation}{0}

\subsection{Total Bohr ergodicity and joint $\epsilon$-period}
\label{sec:2.1}

\pp

Let us first recall the definition of \emph{uniformly almost periodic}
functions and some of their properties (see \cite{Bo}).
Next we will introduce the notions of total Bohr ergodicity and a joint
$\epsilon$-period.

\pp Let $f$ be a real or complex valued function defined on the real
line. A number $\tau$ is called a translation number of $f$
belonging to $\epsilon\geq 0$ (or an $\epsilon$-period) if
\[
     \sup_{x \in \mathbb{R}} |f(x+\tau) - f(x)| \leq \epsilon.
\]
We say that $f$ is a uniformly almost periodic (u.a.p.)
function if it is continuous and if for any $\epsilon >0$ the set
of its translation numbers belonging to $\epsilon$ is
relatively dense (i.e. there exists a number $\ell>0$ such that
any interval of length $\ell$ contains at least one such
translation number). H. Bohr proved that the space of all u.a.p.
functions is a closed sub-algebra of the Banach algebra
$C_b(\mathbb{R})$ of bounded continuous functions equipped with
the uniform norm and that it is the closure of the space of all
(generalized) trigonometric polynomials of the form $$
      \sum_{\rm finite} A_n e^{i \Lambda_n x}
      \qquad (A_n \in \mathbb{C}, \Lambda_n \in \mathbb{R}).
$$
For any u.a.p. function $f$, as is proved by Bohr, the following
limit exists $$ \mathbb{M} f = \lim_{T \rightarrow \infty}
\frac{1}{T} \int_0^T f(x) dx.
 $$
 It is called the {\em Bohr mean} of $f$.
For any locally integrable  not necessarily u.a.p. $f$, we define $\mathbb{M} f$
as the limsup instead of the limit.

\begin{definition}
A  sequence of real numbers $(u_n)_{n \ge 0}$ is said to be {\em
totally Bohr ergodic} if  for any  arithmetic subsequence
$(u_{am+b})_{m\ge 0}$ ($a \ge 1, b\ge 0$ being fixed) and for any real $p>0$,
the sequence $(u_{a m+b}x)_{m\ge 0}$ is uniformly distributed (modulo $p$)
for almost every $x \in \mathbb{R}$ with respect to the Lebesgue measure.
\end{definition}

The following is the main  property of totally Bohr ergodic
sequences that we will use.

\begin{lemma}\label{bohrbeta} Suppose that $(u_n)_{n \ge 0}$ is
a totally Bohr ergodic sequence. Then for any
  u.a.p. function $f$ and any integers $a\geq 1, b\ge 0$, we have
$$
    \lim_{N\rightarrow \infty} \frac{1}{N} \sum_{n=0}^{N}
        f(u_{an+b} x ) = \mathbb{M} f
        \qquad a.e.
$$ \end{lemma}

\begin{proof} It is a consequence of  the fact that $u_{an+b} x$
is uniformly distributed ($\!\!\mod p$) for almost all $x$, for
any real $p>0$, the fact that $f$ can be uniformly approximated
by trigonometric polynomials and the Weyl criterion. \end{proof}

\begin{remark}\label{rem2.1}
{\rm
 Suppose the $(u_n)_{n \ge 0}$ is a sequence such that
  $\inf_{n\ne m}|u_n-u_m|>0$, then the sequence
 $(u_nx)$ is uniformly distributed for almost every point $x$ \cite{KL}.
Consequently, the sequence $(u_n)$ is totally Bohr ergodic.
A more
special case is $u_n = \beta^n$ with $\beta >1$. This is the most
interesting case for us.
On the other hand, no bounded sequence can be totally Bohr ergodic.
}
\end{remark}

\begin{definition}
     Let $(F_n)_{n \ge 0}$ be a sequence of u.a.p. functions.
Let $\epsilon >0$  and $N\in\NN$.
     A real number $\tau$ is called a joint $\epsilon$-translation number for
     $(F_n)_{n\ge N}$ if
     \[
            \sup_{n\ge N} \sup_{x\in \mathbb{R}} |F_n(x +\tau) -
            F_n(x)|\le \epsilon.
     \]
     If for any $\epsilon>0$ there exists $N(\epsilon)$ such that
such joint $\epsilon$-translation numbers for $(F_n)_{n\ge N(\epsilon)}$
are relatively dense, we say that $(F_n)_{n\ge 0}$ has joint periods.
\end{definition}

\subsection{Kingman's theorem}
\pp

 Following ideas of Katznelson and Weiss \cite{KW} we prove the
following version of Kingman's Theorem. The difficulty in our case
is that we have to deal with an infinite measure space.  We
are also dealing with non--stationary sequences.

\begin{theorem}\label{king}
Let  $(u_n)_{n\ge 0}$  be a totally Bohr ergodic sequence of real
numbers and $(f_n)_{n \ge 0}$ be a sequence of uniformly almost
periodic functions. Suppose\\
\indent {\rm (i)} \ The sequence  $(n^{-1}f_n)$ has joint
periods.\\
\indent {\rm (ii)} \  The following  subadditivity is fulfilled
\[
 f_{n+m} (x) \le f_n(x) + f_m(u_nx)
\quad\text{for a.e. $x$ and all $n,m$}.
\]
\indent {\rm (iii)} \ For any $n \ge 1$
\begin{equation} \label{supadd2}
\sup_m \left( f_{m}(u_n x)- f_{m}(x)\right) < \infty \quad\text{for a.e. $x$}.
\end{equation}
Then the following limit exists and is a constant
\[
\lim_{n\to\infty} \frac{1}{n} f_n(x) = \inf_n \frac{1}{n}
\MM f_n
\quad\text{for a.e. $x$}.
\]
\end{theorem}

\begin{proof} \
The proof is a modification of Katznelson-Weiss' proof \cite{KW}.
Without loss of generality we assume that $u_0=1$.

Let $\gamma = \inf_n \frac1n \MM f_n$.
Let us put
\[
f^-(x) = \liminf_{n\to\infty}\frac1n f_n(x), \qquad
 f^+(x) = \limsup_{n\to\infty}\frac1n f_n(x).
 \]
We remark that the subadditivity implies that $f^\pm(x)\le
f^\pm(u_nx)$ for all $n\in\NN$ and a.e. $x\in\RR$. In a finite
measure space this would imply the invariance a.e. In our case
it is the boundedness~\eqref{supadd2} which implies that
$f^\pm(u_n x)\le f^\pm(x)$ a.e.,
what makes the functions $f^\pm$ invariant in the sense that
$f^\pm(x)= f^\pm(u_n x)$ a.e. ($\forall n$).

The first part of the proof, i.e. $f^+\le \gamma$ a.e., is simple.
We just exploit the fact
that any (infinite) arithmetical subsequence of $(u_nx)$ is Bohr-uniform
distributed. This provides us with a kind of ergodic theorem.
In fact,
fix an integer $N$. For any integer $n$ write $n=mN+r$ with $0\le
r<N$. We have
\[
 f_n(x) \le \sum_{k=0}^{m-1}f_N(u_{kN}x) + f_r(u_{mN}x).
\]
The $N$ functions $f_r$ ($r=0,1,\cdots, N-1$) being bounded, by
Lemma~\ref{bohrbeta}, this readily implies by the Bohr-uniform
distribution of $(u_{kN}x)_{k\ge 0}$ that
\[
f^+(x) \le \lim_{m \rightarrow \infty} \frac{1}{mN + r}
\sum_{k=0}^{m-1}f_N(u_{kN}(x)) = \frac1N \MM f_N
\]
for a.e. $x$. Hence $f^+(x)\le \gamma$ for a.e. $x$.

Next, we want to prove that $f^-(x)\ge\gamma$ for a.e. $x$. For
this we assume  that $\gamma>-\infty$, otherwise it is trivially
true. Adding to each $f_n$ the constant value $-n\|f_1\|$ creates
a new subadditive sequence $\tilde f_n:= f_n - n\|f_1\|$ with
$\tilde f_n\le 0$ ($n \ge 1$), $f^-=\tilde f^-+\|f_1\|$ and
$\gamma_f=\gamma_{\tilde f}+\|f_1\|$.
So, we may assume that $f_n\le 0$ ($n \ge 1$).
We can furthermore set $f_1=0$ (observe that this will
not affect the subbadditivity condition since $f_n\le 0$).
Then for any
$\Delta>0$ we truncate the function $f_n$ in the way
\[
 f_{n,\Delta}=\max(f_n,-n\Delta).
\]
Note that the sequence $f_{n,\Delta}$ fulfills the assumptions
of the theorem. Note that in this case $\gamma_\Delta\ge -\Delta$
and also $f^-_\Delta(x)\ge -\Delta$ for all $x$. It is clear
that
\[
      f_\Delta^\pm (x) = \max (f^\pm(x), - \Delta).
\]
Assume that we proved the theorem for the sequence
$f_{n,\Delta}$ for any $\Delta$. Then we claim that
$f_\Delta^-(x)\searrow f^-(x)$ for a.e. $x$ as $\Delta$ goes
to $\infty$. In fact, if $f_\Delta^-(x)>-\Delta$ for some
$\Delta$ then $f_\Delta^-(x)=f^-(x)$
and $\gamma_\Delta = \gamma_f >-\infty$. On the other hand if
$f_\Delta^-(x)=-\Delta$ for all $\Delta$ then obviously
$f^-(x)=-\infty$ and $\gamma_\Delta = \gamma_f =-\infty$.
This proves the theorem for the sequence $f_n$.

>From now on we assume that $f_1=0$, $f_n\le 0$ and  $f_n$
is  truncated and we skip
the subscript $\Delta$. Let $\epsilon>0$. By the hypothesis (i)
on the joint periods, there is an integer $N(\epsilon)$ such that
the joint $\epsilon$-translation numbers are relatively dense.
For these numbers $\tau$ we have
\begin{equation}~\label{E-transl}
               \left| \frac{f_n(x + \tau )}{n} -
               \frac{f_n(x)}{n}\right|
        \le \epsilon \qquad (\forall x \in \mathbb{R}, \forall n \ge N(\epsilon))
\end{equation}
Notice that there is no loss of generality to suppose that $N(\epsilon)$
increases as $\epsilon$ decreases to $0$.
We define
\[
 n_\epsilon(x) = \min\left\{n\ge N(\epsilon)\colon \frac1nf_n(x) \le f^-(x)+\epsilon\right\}
\]
Let $A_K^\epsilon = \{x\colon n_\epsilon(x)>K\}$. Notice that if
$\epsilon'< \epsilon''$, we have $N(\epsilon')\ge N(\epsilon'')$
and $n_{\epsilon'}(x) \ge n_{\epsilon''}(x)$, so that
$A_K^{\epsilon''} \subset A_K^{\epsilon'}$ for $ K >
N(\epsilon')$.

We claim that

\begin{lemma}~\label{meas-A} For any $\epsilon >0$, we have
\[
 \lim_{K\to\infty}\hat\MM (A_K^\epsilon)=0
\]
where $\hat\MM A = \MM 1_A$ denotes the Bohr mean of the
characteristic function of the set $A$
(defined if necessary with the limsup).
\end{lemma}

In order to prove this Lemma~\ref{meas-A} we need the following
lemma which says that $A_K^\epsilon$ is to some extent periodic.

\begin{lemma}~\label{period-A}
   For any joint $\epsilon$-translation  number $\tau$ of
$(f_n/n)_{n \ge N(\epsilon)}$, we have
$A_K^{2 \epsilon} + \tau \subset A_K^{ \epsilon}$.
\end{lemma}

Let us first prove Lemma~\ref{period-A}. Suppose $ x \in
A_K^{2\epsilon} + \tau$, i.e. $x- \tau \in A_K^{2\epsilon}$, then
\[
      \frac{f_n( x - \tau)}{n} > f^-(x-\tau) + 2 \epsilon
      \qquad (N(2\epsilon) \le n \le K).
\]
This, together with the fact that $\tau$ is a joint
$\epsilon$-translation number for all $f_n/n$ with $n \ge
N(\epsilon) (\ge N(2 \epsilon))$ (see~\eqref{E-transl}), implies
\[
     \frac{f_n(x)}{n}
     \ge \frac{f_n(x-\tau)}{n} - \epsilon
     > f^-(x-\tau) +  \epsilon
     \qquad ( N(\epsilon) \le n \le K).
\]
That means $ x \in A_K^\epsilon$. Thus we have finished the proof
of Lemma~\ref{period-A}.

Now let us prove Lemma~\ref{meas-A}. Since joint
$\frac\epsilon2$-translation numbers are relatively dense
there exists $L=L(\frac\epsilon2)>0$ such that any interval of
length $L$ contains such a joint $\frac\epsilon2$-translation number.
Since $\cap_K A_K^{\frac\epsilon2}=\emptyset$,
for any $\eta >0$ there exists $K_0>0$ such that
\begin{equation}~\label{L-ineq}
\left|A_K^{\frac\epsilon2} \bigcap [-L, L]\right| < L\eta \qquad
(\forall K \ge K_0)
\end{equation}
where $|\cdot|$ denotes the Lebesgue measure (see the definition
of $n_\epsilon(x)$). We claim that $\hat\MM(A_K^\epsilon)\le\eta$
($\forall K \ge K_0$). Otherwise $\hat\MM(A_K^\epsilon) > \eta$
for some $K \ge K_0$. Then by the definition of $\hat\MM(A_K^\epsilon)$
there exists $x_0 \in \mathbb{R}$ such that
\[
    \int_{x_0}^{x_0 +L} \chi_{A_K^\epsilon} (x) d x \ge L\eta.
\]
Take a joint $\frac\epsilon2$-translation number $\tau \in [-x_0-L, -x_0
]$, i.e. $-L\le x_0 + \tau \le 0$.  Then by Lemma~\ref{period-A},
we have
\begin{eqnarray*}
    \left|A_K^{\frac\epsilon2} \bigcap [-L, L]\right|
    & \ge & \int_{x_0+\tau}^{x_0+\tau +L} \chi_{A_K^{\frac\epsilon2}}(x) dx \\
       &=& \int_{x_0}^{x_0 +L} \chi_{A_K^{\frac\epsilon2}}(y+\tau) dy\\
       & \ge& \int_{x_0}^{x_0 +L} \chi_{A_K^{\epsilon}}(y) dy
\end{eqnarray*}
For the first inequality we have used the fact that $[x_0+\tau,
x_0+\tau + L] \subset [-L, L ]$ and for the last inequality we
have used Lemma~\ref{period-A}. What we have deduced contradicts
~\eqref{L-ineq}. Thus Lemma~\ref{meas-A} is proved.

We continue our proof of Theorem~\ref{king}.
Let $ S := \|f^-\|_\infty < \infty$.
Let $K$ be such that $\hat\MM(A_K^\epsilon)\le \epsilon/S$.
We define first
\[
 g(x) =
\begin{cases}
f^-(x) & \text{ if }  x\not\in A_K^\epsilon\\
0 & \text{ if }  x\in A_K^\epsilon
\end{cases}
\quad\text{and}\quad
m(x)= \begin{cases}
n^\epsilon(x) & \text{ if }  x\not\in A_K^\epsilon\\
1 & \text{ if }  x\in A_K^\epsilon
\end{cases}
\]
Lemma~\ref{meas-A} implies that
\begin{equation}\label{bm1}
\MM g \le \MM f^- + \epsilon\quad\text{a.e.}.
\end{equation}
Moreover by the invariance of $f^-$ we have (remember that $u_0=1$)
\begin{equation}\label{*2}
g(x) \le g(u_kx)\quad\text{for a.e. $x$ and all $0\le k\le m(x)-1$}.
\end{equation}
Then we have
\begin{equation}\label{*3}
 f_{m(x)}(x) \le (g(x) + \epsilon)m(x)
\le \sum_{k=0}^{m(x)-1} g(u_kx)+\epsilon m(x).
\end{equation}
We define inductively $m_0(x)=0$ and
\[
 m_k(x) = m_{k-1}(x) + m(u_{m_{k-1}(x)}x).
\]
Now choose $R>K$ and let $k(x)$ be the maximal $k$ for which $m_k(x)\le
R$.
Note that $R-m_{k(x)}(x)<K$. Now we get by the subadditivity and
Equation \eqref{*3}
\[
\begin{split}
 f_R(x)
&\le
\sum_{k=0}^{k(x)-1} f_{m(u_{m_k(x)}x)}(u_{m_k(x)}x)
+ \underbrace{f_{R-m_{k(x)}(x)}(u_{m_{k(x)}(x)}x )}_{\le 0}  \\
&\le
\sum_{k=0}^{k(x)-1} \sum_{j=m_{k-1}(x)}^{m_{k}(x)} g(u_jx) +
(m_k(x)-m_{k-1}(x))\epsilon \\
&\le
\sum_{j=0}^{m_{k(x)}(x)} g(u_jx) + m_{k(x)}(x) \epsilon \\
&\le
\sum_{j=0}^{R-1} g(u_jx) + R \epsilon + KS.
\end{split}
\]
Taking the Bohr mean, using that $\MM g=\MM(g\circ u_j)$
and dividing by $R$ gives
\[
\begin{split}
 \frac1R \MM f_R
&\le \MM g + \epsilon + \frac{KS}R \\
&\le \MM f^- +2\epsilon + \frac{KS}R
\end{split}
\]
by Equation \eqref{bm1}. Now we let $R\to\infty$ and we get
\[
 \gamma \le \MM f^-.
\]
We claim that  $f^{-}\le\gamma$  implies
$f^{-}=\gamma$ for a.e. $x$.
Suppose this was not the case, then one could find $\epsilon>0$,
$\delta>0$ and an interval $J=(0,L)$ of length $|J|=L=L_\epsilon$ such
that $|A\cap J|=\delta>0$, where
\[
 A = \{x\colon f^-(x)<\gamma -\e\}
\]
By the invariance of $f^-$ we have $u_kA=A$ for all $k\in\NN$.
Hence, for all $k\in\NN$
\[ \frac{1}{u_kL}\int_0^{u_kL}\chi_{A}\,
 dx=
\frac{1}{L}\int_{0}^{L}\chi_{A}\, dx >\frac{\delta}{L}>0.
\]
Since
$\limsup_k u_k=+\infty$ (see remark \ref{rem2.1}), we have
 $\hat\MM A>\frac\delta L$,
and thus
\[
\underline \MM f^-<\gamma\left(1-\frac{\delta}{L}\right)+(\gamma-\e)\frac{\delta}{L}<\gamma.
\]
\end{proof}

\begin{remark}{\rm
One can prove a similar theorem for more general sequences
 $(u_n(x))$. In this case it seems to be necessary to assume $L^1$
 Bohr-uniform distribution.}
\end{remark}

\subsection{Oseledec theorem}
\setcounter{equation}{0}
 \pp

 Kingman's theorem implies the following Oseledec type theorem
(see Ruelle \cite{Ru}).

\begin{theorem}\label{ose}
Let  $\beta>1$ be a real number.
Let $M\colon \mathbb{R}\to GL_d(\CC)$ be a uniformly almost
periodic function. Write
\[
M_x^n = M(\beta^{n-1}x)\cdots M(\beta x)M(x).
\]
Suppose the $q$-exterior products $\frac{1}{n} \log \|(M_x^n)^{\wedge q}\|$ have
joint periods, for $q=1,\ldots,d$. Then there is $\Gamma\subset
 \mathbb{R}$ with $\beta\Gamma\subset \Gamma$ of full Lebesgue measure
(in the sense that $\RR\setminus\Gamma$ has 0 measure) such
that if $x\in\Gamma$ then

a) $\lim_{n\to\infty} (M_x^{n*}M_x^{n})^{\frac{1}{2n}} = \Lambda_x$ exists.

b) Let $\exp\lambda_x^{(1)}<\cdots<\exp\lambda_x^{(s)}$ be the
 eigenvalues of $\Lambda_x$ (where $s=s(x)$ and the $\lambda_x^{(r)}$
 are reals), and $U_x^{(1)},\ldots,U_x^{(s)}$ the corresponding
 eigenspaces.
Let $m_x^{(r)}=\dim U_x^{(r)}$. We have
$\lambda_{\beta x}^{(r)}=\lambda_{x}^{(r)}$ and $m_{\beta x}^{(r)}=m_x^{(r)}$
and
\[
 \lim_{n\to\infty} \frac1n \log\|M_x^nv\|= \lambda_x^{(r)}
\quad\text{when}\quad
v\in V_x^{(r)}\setminus  V_x^{(r-1)}
\]
for $r=0,\ldots,s$ where $V_x^{(0)}=\{0\}$ and
$V_x^{(r)}=U_x^{(1)}+\cdots+U_x^{(r)}$.

 c) Moreover $V_x^{(r)}$ depends measurably on $x$ and
$M_x V_x^{(r)}=V_{\beta x}^{(r)}$.

d) In addition the functions $\lambda_{x}^{(r)}$ and $m_x^{(r)}$
 are constant a.e.
\end{theorem}

\begin{proof}
This theorem follows in the standard way from  the a.e. convergence  of
\[
\lim_{n\to\infty}\frac{1}{n}\log\|\left(M^n_x\right)^{\wedge q}\|
\]
for $1\le q\le d$, which in the classical case is insured by
Kingman's theorem.
So we only need to check that for the functions
 $f_n^{(q)}(x)=\log\|\left(M^n_x\right)^{\wedge q}\|$
and the sequence $u_n=\beta^n$ the assumptions
 of Theorem~\ref{king} hold.

First we note that $M^{-1}$ is uniformly almost periodic because
$M$ is uniformly almost periodic and $M^{-1}(x)\in GL_d(\CC)$.
Second we note that $M^{\wedge q}$ and $(M^{\wedge q})^{-1}$ are
again uniformly almost periodic, since each entry is a rational function
of the entries of $M$ and $M^{-1}$, respectively.
Hence,
 \[
 \sup_{x\in\RR}\|M(x)^{\wedge q}\|+\sup_{x\in\RR}\|(M^{\wedge q})^{-1}(x)\|=W_q<\infty.
 \]
  Subadditivity is obviously fulfilled since
$(M^{n+m}_x)^{\wedge q} = (M_{\beta^nx}^m)^{\wedge q}(M_x^n)^{\wedge q}$.

 Condition~\eqref{supadd2} follows from
\[
 \|(M^{n}_{\beta x})^{\wedge q}\|=
\|M(\beta^n x)^{\wedge q} (M^n_x)^{\wedge q} (M(x)^{\wedge q})^{-1}\|
\le \|(M^n_x)^{\wedge q}\|+W_q.
\]
Finally by Remark~\ref{rem2.1} the sequence $\beta^n$ is totally Bohr ergodic,
so Theorem~\ref{king} applies. Assertions (a), (b) and (c) follows from
Proposition 1.3 (see also the proof of Theorem 1.6) in \cite{Ru}.

Now we prove d).
By Kingman's theorem (Theorem \ref{king}), we have
for almost all $x$
\begin{equation}\label{Lyap1}
     \lim_{ n \to \infty} \frac{1}{n} \log \| M_x^ {n \wedge k}\|
         = \inf_{ n \ge 1} \frac{1}{n} \mathbb{M} \log \|M_x^ {n \wedge k}\|
         \qquad (1 \le k\le s ).
\end{equation}
On the other hand, by the properties of exterior product,
if we write $k = \sum_{i=1}^ {j-1} m_x^{(s - i)} + \ell$
with $0\le j<s$ and $ 0\le \ell \le m_x^ {(s-j)}$ we have
\begin{equation}\label{Lyap2}
     \lim_{ n \to \infty} \frac{1}{n} \log \| M_x^ {n \wedge k}\|
         = m_x^ {(s)} \lambda_x^ {(s)} + \cdots +
            m_x^ {(s-j+1)} \lambda_x^ {(s-j +1)} +
                     \ell  \lambda_x^ {(s-j)}.
\end{equation}
We can solve  $\lambda_x^ {(r)}$
($1\le r \le s $) from the system (\ref{Lyap1})-(\ref{Lyap2}).
The solution is independent of $x$ since it depends
 only on the right hand side terms in
(\ref{Lyap1}). Consequently  $m_x^{(r)}$ is also independent of $x$.
\end{proof}

%%%%%%%%%%%%%%%%%%%%%%%%%%%%%%%%%%%%%%%%%%%%%%%%%%%%%%%%%%%%%
\section{When $\beta$ is a Pisot--Vijayaraghavan number}\label{sec:3}

We restrict our attention to the special case where $\beta>1$ is a
Pisot--Vijayaraghavan (PV) number and $f_n$ are defined
by~\eqref{log-product-matrix}.
We will prove that, under some extra condition, the sequence
$n^{-1} f_n$ has joint periods and
the Kingman theorem and the Oseledec theorem apply. To do this, we
need a distortion lemma and some properties of PV-numbers.

\subsection{Distortion lemmas}
\pp

\begin{lemma}\label{dist}
Let $M\colon \RR\to GL_d(\CC)$ such  that
\begin{equation}~\label{def-D}
    D: =\sup_{x\in \mathbb{R}} \|M(x)\| \ \|M(x)^{-1}\| <\infty
    \end{equation}
 Let $(x_k)$ and $(y_k)$ be two sequences in $\mathbb{R}$ and let
$\theta_k=\|M(x_k)-M(y_k)\|$. Then
for all $n\in\NN$ and $0\neq {\bf v} \in \mathbb{C}^d$ we have
\[
 \frac{|M(x_1)\cdots M(x_{n-1}) M(x_n){\bf v}|}
      {|M(y_1)\cdots M(y_{n-1}) M(y_n){\bf v}|}
      \le
 \exp \left( C \sum_{k=1}^n D^k \theta_k \right).
\]
where $C = (\sup_{x \in \mathbb{R}} \|M(x)\|)^ {-1}$.
\end{lemma}

\begin{proof}
Let
\begin{eqnarray*}
    Q_1(x) & = &  \II, \qquad Q_k(x)  = M(x_1)M(x_2)\cdots M(x_{k-1}) \qquad  (1< k\le n)
\\
    Q^n(x) & = & \II, \qquad  Q^k(x)  = M(x_{k+1})M(x_{k+2})\cdots M(x_n) \quad
    (1\le k<n).
\end{eqnarray*}
 We can write
\begin{equation}~\label{quotient}
\frac{ |Q_n(x){\bf v}| } {|Q_n(y){\bf v}|}
= \prod_{k=1}^n
         \frac{ |Q_{k}(y) M(x_k)  Q^k(x){\bf v}|}
              { |Q_{k}(y) M(y_k)  Q^k(x){\bf v}|}.
\end{equation}
Setting $E_k=M(x_k)\cdot M^{-1}(y_k)-\II$
 we get the following estimate for the numerator of the general term in the above
 product
\begin{eqnarray*}
|Q_{k}(y)M(x_k) Q^k(x){\bf v}|
  & = &|Q_k(x) (\I + E_k) M(y_k)
    Q^k(x){\bf v}|\\
&= & |Q_k(x) (\I + E_k) Q_k(x)^{-1} Q_k(x) M(y_k)
    Q^k(x){\bf v}|\\
&= &|(\I + \tilde{E}_k) Q_k(x) M(y_k)
    Q^k(x){\bf v}|
\end{eqnarray*}
where \(
 \tilde {E}_k=Q_k(x)E_k Q_k(x)^{-1}.
\) It follows that
\begin{equation}~\label{quotient2}
\frac{ |Q_{k}(y) M(x_k)  Q^k(x){\bf v}|}
              { |Q_{k}(y) M(y_k)  Q^k(x){\bf v}|}
 \le \| \I +  \tilde{E}_k\|.
\end{equation}
It is obvious that
\begin{equation}~\label{E-E}
 \|\tilde E_k\|\le D^{k-1}\|E_k\|.
\end{equation}
On the other hand
\begin{equation}~\label{Hold}
 \|E_k\| \le \sup_x \|M(x)^{-1}\| \theta_k.
\end{equation}
By combining~\eqref{quotient}, \eqref{quotient2}, \eqref{E-E} and
\eqref{Hold}, we obtain
\[
\frac{|M(x_1)\cdots M(x_{n-1}) M(x_n){\bf v}|}
      {|M(y_1)\cdots M(y_{n-1}) M(y_n){\bf v}|}
 \le\prod_{k=0}^{n-1}\left(1+ C    D^k \theta_k \right).
\]
\end{proof}

If $M(x)$ is non-negative, the next lemma shows that
condition~\eqref{def-D} is not needed for positive vectors

\begin{lemma}\label{dist2}
Let $M\colon \RR\to GL_d(\RR)$ be such that the entries of $M(x)$
are either identically zero or bounded from below by a positive
number $\delta>0$ (independent of entries).
  Then
  for any sequences $(x_k)$ and $(y_k)$ in $\mathbb{R}$ and
  for any non-negative vector ${\bf v}$ we have
\[
 \frac{|M(x_1)\cdots M(x_{n-1}) M(x_n){\bf v}|}
      {|M(y_1)\cdots M(y_{n-1}) M(y_n){\bf v}|} \le
 \exp \left( \frac{1}{\delta} \sum_{k=1}^n  \theta_k \right),
\]
where $\theta_k=\|M(x_k)-M(y_k)\|$
and the norm $|v|= \sum_{i=1}^d|v_i|$ on $\mathbb{R}^ d$ is chosen.
\end{lemma}

\begin{proof}
We may write
\begin{eqnarray*}
   & & |M(x_1)\cdots M(x_{n-1}) M(x_n){\bf v}|\\
    &= & \sum_{i_0, i_1\cdots, i_n} M(x_1)_{i_0, i_1} M(x_2)_{i_1, i_2}
       \cdots M(x_n)_{i_{n-1}, i_n} v_{i_n}
\end{eqnarray*}
We have a similar expression for $|M(y_1)\cdots M(y_{n-1})
M(y_n){\bf v}|$. Now compare the two expressions term by term. By the
hypothesis, both quantities $M(x_1)_{i_0, i_1}$ and $M(y_1)_{i_0,
i_1}$ are either zero or larger than $\delta$. So, using the
trivial inequality $x/y \le e^{x/y-1}$ we have
\[
            M(x_1)_{i_0, i_1} \le M(y_1)_{i_0, i_1} e^{\delta^{-1}\theta_k }.
\]
The same estimates hold for other pairs
$M(x_k)_{i_{k-1}, i_k}$ and $M(y_k)_{i_{k-1}, i_k}$.
The desired inequality follows.
\end{proof}

\subsection{Two properties of PV--numbers}
\pp

 Let $\beta>1$ be a PV--number of order $r$. We  denote its
conjugates by $\beta'_1, \cdots, \beta'_{r-1}$.  Then for $n\geq
1$, denote $$ F_n = \beta^n + {\beta'}_1^n + \cdots +
{\beta'}_{r-1}^n. $$

\begin{lemma}~\label{pisot1}
The number $F_n$ is an integer and we have
       $$
              |\beta^n - F_n| \leq (r-1) \rho^n  \qquad(\forall  n\geq 1)
 $$ where $
            \rho = \max_{1\leq j\leq r-1} |\beta'_j| <1.
$
\end{lemma}

\pp Given any real number $\beta >1$ (not necessarily integral), we can expand
each number $x\in [0,1)$ in a canonical way into its {\em
$\beta$-expansion} \cite{Re} (see also \cite{P,Bl}): $$
      x= \sum_{n=1}^\infty \frac{\epsilon_n(x)}{\beta^n}
$$
 where $(\epsilon_n(x))_{n\geq 1}$ is a uniquely determined sequence in
 $\{0,1,\cdots, [\beta]\}^\mathbb{N}$. We may also call  $(\epsilon_n(x))_{n\geq 1}$
 the $\beta$-expansion of $x$. We note that not  all sequences in  $\{0,1,\cdots, [\beta]\}^\mathbb{N}$
 are $\beta$-expansions. Let $D_\beta$ be the set of all possible
 $\beta$-expansions of numbers in $[0,1)$.   A finite sequence
 $\epsilon=(\epsilon_1, \cdots, \epsilon_n)$
 (of length $n$) in $\{0,1,\cdots, [\beta]\}^n$ is said to be {\em
 admissible}
 if it is the prefix  of the $\beta$-expansion of some number $x$. For
 such an admissible sequence, we define
 $$
     I(\epsilon_1, \cdots, \epsilon_n)
      = \{x\in [0,1):  \epsilon_1(x) = \epsilon_1, \cdots, \epsilon_n(x) =
      \epsilon_n\}.
 $$
 It is known that if $D_\beta$ is endowed with
 the lexicographical order, the map which associates  $x$ to its
 $\beta$-expansion is strictly increasing. The set $I(\epsilon)=I(\epsilon_1, \cdots, \epsilon_n)$
 is an interval, called a {\em $\beta$-interval} of level $n$. Its length
 is denoted by $|I(\epsilon)|$.

\begin{lemma}~\label{pisot2} Suppose $\beta>1$ is a PV--number.
 There is a constant $C >0$ such that
$$
       C^{-1} \beta^{-n}  \leq   |  I(\epsilon_1,\cdots,\epsilon_n) |  \leq C
       \beta^{-n}
$$
 for any integer $n\geq 1$ and any $\beta$-interval
 $ I(\epsilon_1,\cdots,\epsilon_n)$.
\end{lemma}

See \cite{F} for proofs of Lemma \ref{pisot1} and Lemma \ref{pisot2}.

\subsection{Existence of joint periods}
\pp

\begin{definition} Let $\beta >1$ be a positive real number and let $M:
\mathbb{R} \rightarrow \mathcal{M}_{d \times d}(\mathbb{C})$.   If
the entries of $M$ are functions of the form $f(\beta^n x)$ ($n
\in \mathbb{Z}$) where $f$ is $1$--periodic continuous, we say that
$M$ is {\em $\beta$--adapted u.a.p.}
\end{definition}

\begin{remark}{\rm The matrix $M(x)$ defined by
\eqref{Matrix-multip} associated to a multiperiodic function
is $\beta$--adapted. }
\end{remark}

\begin{proposition}~\label{J-period}
      Let $\beta >1$ be a PV--number and let
      $M : \mathbb{R} \rightarrow \mbox{GL}_d(\mathbb{C})$
      be $\beta$--adapted and $\alpha$--H\"{o}lder continuous.
      Suppose that
       \begin{equation}~\label{cond-D-rho}
              D\rho^\alpha <1
\end{equation}
where $\rho$ is the maximal modulus of the conjugates of $\beta$ and $D$
is the same as in the distortion lemma (Lemma~\ref{dist}).
Then for any $1\le q\le d$ the sequence $n^{-1} f_n^{(q)}(x)$ has joint
periods, where
$$
f_n^{(q)} (x) = \log \|(M(\beta^{n-1} x) \cdots
M(\beta x) M(x)))^{\wedge q}\|.
$$
\end{proposition}

\begin{proof}
Since $M$ is $\beta$--adapted,  the entries of $M(\beta^k x)$ are all of the
form $h_{i,j}(\beta^{\ell_{i,j}} x)$ with 1--periodic function $h_{i,j}(x)$
and integer $\ell_{i,j}\ge 0$ for sufficiently large $k$.
So, if necessary, we  consider
$\log \|(M(\beta^{n-1} x)\cdots M(\beta^{k_0} x))^{\wedge q}\|$
for some sufficiently large but fixed $k_0\ge 0$.

Consider $\tau = \beta^m \eta_m + \cdots + \beta \eta_1 + \eta_0$
where $m \ge 1$ and $0\le \eta_i \le \beta$ are integers. We are
going to show that all such $\tau$ are joint
$\epsilon$-translation numbers for $n^{-1} f_n^ {(n)}(x)$ with $n \ge
N(\epsilon)$, where $N(\epsilon)$ depending on $\epsilon$ is an
integer to be determined.

By Lemma~\ref{pisot1}, we have
\[
       \inf_{j \in \mathbb{Z}} |\beta^k \tau - j| \le C'\rho^k
\]
for all $k$ and some constant $C'$ independent of $k$ and $\tau$.
For $k\ge 0$ each entry of $M^{\wedge q}(\beta^{k+k_0}x)$
is a degree $q$ polynomial in $d^2$ variables of the form
$h(\beta^{\ell+k}x)$, with $h$ $\alpha$-H\"older, 1--periodic, and $\ell\ge0$.
Notice that we have
$h(\beta^{\ell+k}(x+\tau))=h(\beta^{\ell+k}x)+O(\rho^{k\alpha})$,
hence
\[
\|M(\beta^{k+k_0}(x+\tau))^{\wedge q}-M(\beta^{k+k_0}(x))^{\wedge q}\|
= C_q \rho^{k\alpha}
\]
for some constant $C_q$.
By the distortion Lemma~\ref{dist} and the above estimate we have
\begin{eqnarray*}
    & &|f_n^{(q)}(x+\tau) - f_n^{(q)}(x)|\\
    & = & \left|\log
    \frac{\| (M(\beta^{n-1} x + \beta^{n-1} \tau)\cdots M(\beta x + \beta \tau) M(x + \tau))^{\wedge q}\|}
        { \|(M(\beta^{n-1} x) \cdots M(\beta x) M(x))^{\wedge q}\|}\right|\\
    & = &  C_q C'\sum_{k=1}^n D^k \rho^{(k-k_0)\alpha}
\le \frac{CC'D\rho^{-k_0}}{1- D \rho^\alpha}
    =: \mathcal{C}
\end{eqnarray*}
So, we may choose $N(\epsilon)= \mathcal{C}/\epsilon$. In order to
finish the proof, it suffices to notice that Lemma~\ref{pisot2}
implies that all these $\tau$ form a subset with bounded gap in
$\mathbb{R}$.
\end{proof}

\begin{proposition}~\label{J-period-positive}
      Let $\beta >1$ be a PV--number and let
      $M : \mathbb{R} \rightarrow \mbox{GL}_d(\mathbb{C})$
      be $\beta$--adapted and $\alpha$--H\"{o}lder continuous.
      Suppose that the entries of $M(x)$ are either identically
      zero or larger than a constant $\delta>0$.
      Then $n^{-1} f_n(x)$ has
joint periods, where $$f_n^{(q)} (x) =
\log \|( M(\beta^{n-1} x) \cdots
M(\beta x)M(x))^{\wedge q})\|.$$
\end{proposition}

\begin{proof}
   The proof is the same as the last proposition. But we use
   Lemma~\ref{dist2} instead of Lemma~\ref{dist}.
\end{proof}

\bigskip

\bigskip
\section{Multiperiodic functions}
\setcounter{equation}{0}

 \pp As we pointed out in the introduction
and as we will see in Section~\ref{sec4.2}, our scalar
equation~\eqref{eq-multiperiodic} can be converted to the vector
equation~\eqref{eq-vector}. So, we first study the vector
equation~\eqref{eq-vector}.

\subsection{Equation $G(x) = M(x/\beta)G(x/\beta)$} \pp

 Let $M: \mathbb{R} \rightarrow \mathcal{M}_{d \times d}(\mathbb{C})$
be a  matrix valued function. We consider  the following vector
valued
 equation $$
      G(x) = M\left(\frac{x}{\beta}\right)G\left(\frac{x}{\beta}\right)
$$ where the unknown $G: \mathbb{R} \rightarrow \mathbb{C}^d$ is a vector
valued function.

\begin{theorem}\label{existence-G}
Let $\beta >1$ be a real number and $M(x)$ be a complex matrix valued
function. Suppose that $M$ is Lipschitzian and that $M(0)$ is non
negative and has $1$ as a simple eigenvalue with a corresponding
strictly positive eigenvector ${\bf v}$. Then there exists,
up to a multiplicative constant, a
unique continuous solution $G(0)\not=0$ of the equation $G(x) =
M(x/\beta)G(x/\beta)$. The
solution can be defined by
 $$
     G(x) = \lim_{n \rightarrow \infty} M\left(\frac{x}{\beta} \right)
     M\left(\frac{x}{\beta^2} \right) \cdots M\left(\frac{x}{\beta^n} \right)
     {\bf v}
$$ where the convergence  is uniform on every compact subset in
$\mathbb{R}$.
\end{theorem}

\begin{proof}  \  \
Write ${\bf v} = (v_1, \cdots, v_d)^t$. We introduce the following
norm for $\mathbb{C}^d$ $$
               \|z\| = \max_{1\leq j \leq d} \frac{ |z_j|}{v_j}
               \qquad ( z = (z_j)_{1\leq j\leq d} \in \mathbb{C}^d).
$$ Then a matrix $A= (a_{i, j}) \in \mathcal{M}_{d \times d}(\mathbb{C})$,
considered as  an operator on the normed space $(\mathbb{C}^d,
\|\cdot\|)$, admits its operator norm $$ \|A\| = \max_{1\leq i\leq
d}
 \frac{1}{v_i}\sum_{j=1}^d  |a_{i,j}| v_j. $$
 Notice that $\|M(0)\|=1$ because $M(0) {\bf v} = {\bf v}$.

\pp Since the eigenvalue $1$ of  $M(0)$ is simple (and isolated),
and $M(x)$ is Lipschitz continuous, by the perturbation
theory of matrices, there is a neighborhood of $0$, say $[-\delta,
\delta]$ ($\delta >0$), such that for any $x \in [-\delta,
\delta]$,   $M(x)$ has a simple eigenvalue $\lambda(x)$ and a
corresponding  eigenvector $v(x)$ satisfying
\begin{equation}~\label{perturb-eigen}
       |\lambda(x)-1| \leq C|x|,
       \qquad \|v(x) - {\bf v}\| \leq C|x|
       \qquad (x \in [-\delta, \delta])
\end{equation}
 for some constant $C>0$. We claim that
  the  limit
 \begin{equation}~\label{lim-G}
       G(x) = \lim_{n\rightarrow \infty}
        M\left(\frac{x}{\beta}\right) M\left(\frac{x}{\beta^2}\right) \cdots
       M\left(\frac{x}{\beta^n}\right){\bf v}
 \end{equation}
 exists (uniformly on any compact set). It is clear  that
 the limit function is  a solution.

   \pp
   Denote
 $$
     Q_n(x) = M\left(\frac{x}{\beta}\right) M\left(\frac{x}{\beta^2}\right) \cdots
       M\left(\frac{x}{\beta^n}\right).
 $$
 The proof of the existence of the limit in \eqref{lim-G} is based on the
following lemma.

\begin{lemma} For any $\delta>0$, there exists a constant $D>0$ such that
for any $ n \geq 1$ and any $ x\in [-\delta, \delta]$ we have
 $$
        \|Q_n(x) \| \leq D  \qquad
        \|Q_n(x){\bf v} - {\bf v} \| \leq D |x|.
 $$
 \end{lemma}

\pp
 To get the boundedness of $\|Q_n(x) \|$, it suffices to notice that
 $$
         \|Q_n(x) \| \leq \prod_{j=1}^n f\left( \frac{x}{\beta^j}\right)
 $$
 where the scalar function $f(x)= \|M(x)\|$ is Lipschitzian and $f(0)=1$
 (we have used
 our choice of the norm of $\mathbb{C}^d$), and that the products converge uniformly
 on $[-\delta, \delta]$ to a continuous function \cite{FL}.
 Now we  prove that
 \begin{equation}~\label{pre-Cauchy}
    \|Q_n(x){\bf v}
         - Q_{n-1}(x) {\bf v}\| \leq C' \frac{|x|}{\beta^n}
\end{equation}
 where $C'>0$ is some constant.
 In fact, since $M(x) v(x) = \lambda(x) v(x)$, we have
 $$
    M\left(\frac{x}{\beta^n}\right){\bf v} - {\bf v} =  M\left(\frac{x}{\beta^n}\right)
    \left[{\bf v} - v\left(\frac{x}{\beta^n}\right)\right] + \left[
       \lambda\left(\frac{x}{\beta^n}\right) v\left(\frac{x}{\beta^n}\right) - {\bf
    v}\right].
 $$
 Multiplying both  sides by $Q_{n-1}(x)$,
  we get
 \begin{eqnarray*}
       &     & \| Q_n(x){\bf v} -   Q_{n-1}(x){\bf v}\|\\
      & \leq & \left\|Q_n(x) \left({\bf v} - v\left(\frac{x}{\beta^n}\right)\right)\right\|
         + \left\|Q_{n-1}(x) \left({\bf v} -
        \lambda\left(\frac{x}{\beta^n}\right)  v\left(\frac{x}{\beta^n}\right)\right)\right\|
 \end{eqnarray*}
 Notice that
$$
  \|\lambda(x) v(x) - {\bf v}\|
  \leq \|\lambda(x) -1\| \ \|v(x)\| + \|v(x) - {\bf v}\|.
$$ Using the last inequality, the estimates in \eqref{perturb-eigen}
and that we have just proved $\|Q_n(x)\|\leq D$,
we obtain~\eqref{pre-Cauchy}. Then for $n>m$
 $$
       \| Q_n(x){\bf v} -  Q_m(x) {\bf v} \|
          \leq  \sum_{k=m+1}^n \| Q_k(x){\bf v} -   Q_{k-1}(x){\bf v}\|
            \leq
            \frac{C'|x|}{\beta^{m-1}(\beta -1)}.
 $$
 That means $Q_n(x) {\bf v}$ is a Cauchy sequence in the space
  $C([-\delta, \delta])$
 of
 continuous functions equipped with uniform norm.
Since for any fixed integer $n_0$, we have $$
   \lim_{n\rightarrow\infty} Q_n(x) {\bf v} =
      Q_{n_0}(x) \cdot \lim_{n\rightarrow \infty} Q_n\left(\frac{x}{\beta^{n_0}}\right) {\bf
      v},
$$
it follows that the uniform convergence of $Q_n(x) {\bf v}$ on
$[-\delta, \delta]$ implies its uniform convergence on any compact
set.

The uniqueness of solution is easy. Let $G\neq 0$ be a solution.
First notice that $G(0)$ is an eigenvector of $M(0)$ associated
to the simple eigenvalue 1.
Hence we may assume that $G(0)= {\bf v}$.  By iterating the
equation, we get
\begin{eqnarray*}
    G(x)  =  Q_n(x)
    G\left(\frac{x}{\beta^n}\right)
     =  Q_n(x){\bf v} + Q_n(x) \left( G\left(\frac{x}{\beta^n}\right)- {\bf v}\right)
\end{eqnarray*}
The last term converges to zero (uniformly on any compact set)
because of $\|Q_n(x)\| \leq D$. So, $G(x)$ must be the limit of
$Q_n(x){\bf v }$.
 \end{proof}
\medskip

 \pp
 \begin{remark} {\rm
 In the theorem, neither the almost periodicity of $M(x)$ nor the positivity of $M(x)$
 is  required, but only the positivity of $M(0)$. That $1$ is
 an eigenvalue of $M(0)$ is necessary for the equation \eqref{eq-vector}
to have a solution $G(x)$ such that $G(0)\not=0$.
 }
 \end{remark}

 \begin{remark} {\rm
 The Lipschitz continuity is not really
necessary. H\"{o}lder continuity or even Dini continuity is
sufficient. }
\end{remark}

\begin{remark} {\rm
 If the entries of $M$ are (real) analytic, then the
solution $G$ is also analytic. Because, for any $x_0 \in
\mathbb{R}$, there is a disk on the complex plane centered at $x_0$
on which  $Q_n(x) {\bf v}$ (as functions of complex variable $x$)
uniformly converges. }
\end{remark}

\subsection{Existence of multiperiodic functions}\label{sec4.2}
\pp

Here we give a proof of Theorem A based on
Theorem~\ref{existence-G}.

Let $M(x)$ be as in \eqref{Matrix-multip}.
It is easy to see that the characteristic polynomial of $M(0)$
takes the form $$
    P(u) = u^d - f_1(0) u^{d-1} - f_2(0) u^{d-2} - \cdots - f_{d-1}(0) u - f_d(0)
$$ The consistency condition implies that $1$ is an eigenvalue of $M(0)$.
Notice that $$
  P'(1) = f_1(0) +2 f_2(0)  \cdots + (d-1) f_{d-1}(0)  + d f_d(0) >0.
$$ So, the eigenvalue $1$ is simple. By Theorem~\ref{existence-G}, there is a unique
solution of $G(x) = M(x/\beta)G(x/\beta)$. Let
 $$G(x) = \left(
                               G_1(x), G_2(x), \cdots, G_d(x)
\right)^t $$ Then $G_1(x)$ is a solution of \eqref{eq-multiperiodic}.
 If $F$ is a solution of \eqref{eq-multiperiodic}. Let $$
   \tilde{G}_1(x) = F(x),   \tilde{G}_2(x) = F(x/\beta),\cdots,  \tilde{G}_d(x)=
     F(x/\beta^{d-1}).
 $$
 Then $ \tilde{G}= ( \tilde{G}_1, \cdots,  \tilde{G}_d)^t$ is a solution of
 $G(x) = M(x/\beta) G(x/\beta)$. Thus the uniqueness of the solution of
Equation~\eqref{eq-vector} implies that of
Equation~\eqref{eq-multiperiodic}.

\subsection{Asymptotic behavior  of multiperiodic functions}\label{sec4.3}
\pp

Let us consider the asymptotic behavior of a multiperiodic
function, or more generally the asymptotic behavior of a
solution $G$ of Equation~\eqref{eq-vector} provided it exists
(the existence may be guaranteed by Theorem~\ref{existence-G}.

\medskip
\begin{theorem}~\label{asymp-pointwise}
Let $\beta >1$ be a PV--number
whose maximal conjugate has
modulus $\rho$. Let $M: \mathbb{R} \to {\rm GL}_d(\mathbb{C})$
be a $\beta$--adapted u.a.p. H\"{o}lder function of order $\alpha>0$.
Suppose that $G$ is a solution of $G(x)= M(x/\beta) G(x/\beta)$.
Suppose furthermore that one of the following conditions is satisfied

(i) $D \rho^\alpha <1$ where
$D =\sup_{x\in \mathbb{R}} \|M(x)\| \|M(x)^{-1}\|$
(NB. $\beta$ must be Pisot).

(ii) The entries of $M(x)$ are either identically zero or larger than
a constant $\delta>0$.

Then for a.e. $x\in\RR$ the limit
\[
h(x)=\lim_{n\to\infty}\frac1n\log |G(\beta^nx)|
\]
exists and  is independent of $x$.
\end{theorem}

\begin{proof}
We first consider the case (i).
By Proposition~\ref{J-period}, Theorem~\ref{ose} applies.
Hence for a.e. $x$, if we denote by $r(x)$ the integer such
that the vector $G(x) \in V_x^{(r)}\setminus V_x^{(r-1)}$ we get
\[
\lim_{n\to\infty} \frac 1n \log |G(\beta^nx)| =
\lim_{n\to\infty} \frac 1n \log |M_x^n G(x)| =
\lambda_x^{(r(x))}.
\]
But $G(\beta x)=M(x)G(x)$ hence
$G(\beta x)\in V_{\beta x}^{(r)}\setminus V_{\beta x}^{(r-1)}$,
from what follows $r(\beta x)=r(x)$, i.e. $r$ is invariant.
Hence constant a.e. because of the total Bohr ergodicity of the
sequence $\beta^n$.
\\

Case (ii).
We use the notation of Proposition~\ref{J-period-positive}.
Since $n^{-1}f_n$ has joint periods, Theorem~\ref{king} applies
(see the proof of Theorem~\ref{ose} for details). Hence the following
limit exists a.e.
\[
\lim_{n\to\infty} \frac 1n f_n(x) = \mathcal{L},
\]
where $\mathcal{L} = \inf_n \frac 1n \MM(f_n)$.
In view of $G(\beta^nx) = M(\beta^{n-1}x)\cdots M(x) G(x)$ the
positivity of $M$ and $G$ gives
\[
\lim_{n\to\infty}\frac 1n\log|G(\beta^nx)| = \mathcal{L} \quad\text{a.e.}
\]
\end{proof}

Note that when $G$ is the solution of equation \eqref{eq-vector} with
$M$ and $G$ given by \eqref{Matrix-multip} we have
\begin{equation}\label{GF}
       |G(x)| = |F(x)| + |F(x/\beta)| + \cdots +
       |F(x/\beta^{d-1})|,
\end{equation}
thus the asymptotic behavior of $\frac 1n \log|G(\beta^nx)|$ and
$\frac 1n \log \sum_{j=0}^{d-1} |F(\beta^{n-j}x)|$ are the same.
Thus Theorem~C follows as an immediate corollary of
Theorem~\ref{asymp-pointwise}.
This partially answers  a question in \cite{JRS} (Conjecture 4.1., p. 263).

We prove now Theorem B. By the primitivity of $M(0)$
and the hypothesis, there exists an ingeter $\tau \ge 1$ such that
$\tilde{M}(x) := M(x/\beta^{\tau-1}) \cdots M(x/\beta)M(x)$ has all
its entries strictly positive (even larger than $c \delta^\tau$ for
some constant $c>0$). Consider the equation
$$
    G(x) = \tilde{M}(x/\beta) G(x/\beta^\tau).
$$
We examine the first entries of both sides. We can find two constants
$0<c_1 <c_2$
such that
we get
$$
   c_1 F(x) \le  F(x/\beta^{\tau +1}) + \cdots +  F(x/\beta^{\tau +d})
       \le c_2 F(x).
$$
Thus Theorem B follows from Theorem \ref{asymp-pointwise}.

\medskip
\begin{theorem}~\label{asymp-mean}
Under the same conditions as Theorem~\ref{asymp-pointwise}, for any
$q \in \mathbb{R}^+$,
 the following limit exists
 $$
      \lim_{n \rightarrow \infty} \frac{1}{n} \log \int_0^1
              \|M(\beta^{n-1} x) \cdots M(\beta x) M(x)\|^q d x.
 $$
\end{theorem}

\begin{proof}
 Write
 $$
 Z_n = \int_0^1 P_n(x)^q  dx \quad {\rm with}\quad
      P_n(x)= \|M(\beta^{n-1} x) \cdots M(\beta x) M(x)\|.
$$
 It suffices to show that there is a
constant $C>0$ such that
 $$
      Z_{n+m} \leq C Z_n Z_m \qquad (n \geq 1, m \geq 1).
 $$
 We assume that $q=1$, just for simplicity. We will use the fact that
 there is a constant $L>0$ such that $$
    \|M(x) - M(y)\|_2 \leq L |x -y|^\alpha
    \qquad (\forall x, y \in \mathbb{R}).
$$
We use the notation
$\Pi_{i=0}^{n} M_i = M_n M_{n-1} \cdots M_1 M_0$
for the (noncommutative) product of the matrices $M_0,\ldots, M_n$.
Write
\begin{eqnarray*}
    Z_{n+m} & =& \sum_{\epsilon}
                    \int_{I(\epsilon)}
                    \left\|
 \prod_{k=0}^{m-1}  M(\beta^{n+k} x)
\cdot  \prod_{j=0}^{n-1}  M(\beta^j x)  \right\|
 dx\\
      & \leq &
          \sum_{\epsilon}
                    \int_{I(\epsilon)}
                   \left\| \prod_{j=0}^{n-1}  M(\beta^j x)\right\| \cdot
                           \left\| \prod_{k=0}^{m-1}  M(\beta^{n+k} x)\right\| dx
\end{eqnarray*}
 where the sum is taken over all $\beta$-intervals
$I(\epsilon)$ of level $n$ (see Lemma~\ref{pisot2}).
Let $a_\epsilon$ be the left endpoint of $I(\epsilon)$.
The integral in the last sum , after the change
of variables $\beta^n (x-a_\epsilon) = y$, becomes $$
   \beta^{-n} \int_0^{\beta^n |I(\epsilon)|}
            \left\| \prod_{j=0}^{n-1}  M(\beta^j a_\epsilon +\beta^{-n +j }y) \right\| \cdot
                          \left\|   \prod_{k=0}^{m-1}  M(\beta^ k y + \beta^{n+k} a_\epsilon)
                          \right\| dy.
 $$
 Notice that if $\epsilon = (\epsilon_1, \cdots, \epsilon_n)$, then
$$
    \beta^{n+k}a_\epsilon =  \beta^{n+k} \left(\frac{\epsilon_1}{\beta} +
    \cdots + \frac{\epsilon_n}{\beta^n}\right)
    = \beta^{n+k-1} \epsilon_1 + \cdots + \beta^k \epsilon_n.
$$
So, by Lemma~\ref{pisot1}, there is an integer $n_\epsilon$ such that
\begin{equation}\label{2.4}
      |\beta^{n+k} a_\epsilon - n_\epsilon|= O
      (\rho^k + \rho^{k+1} + \cdots + \rho^{n+k-1})
        = O(\rho^k)
\end{equation}
By the distortion lemma (Lemma~\ref{dist}, Lemma~\ref{dist2}), we have
 $$
 \left\| \prod_{j=0}^{n-1}  M(\beta^j a_\epsilon +\beta^{-n +j }y)
 \right\|
  \leq C  \left\| \prod_{j=0}^{n-1}  M(\beta^j a_\epsilon) \right\|
$$
$$ \left\|   \prod_{k=0}^{m-1}  M(\beta^ k y + \beta^{n+k} a_\epsilon)
                       \right\|
                       \leq C
                       \left\|   \prod_{k=0}^{m-1}  M(\beta^ k y)
                          \right\|
$$
Therefore, we get
$$
    Z_{n+m}  \leq C  \beta^{-n}\sum_{\epsilon}  P_n  ( a_\epsilon)
                    \leq C'Z_n Z_m.
$$
 \end{proof}

\begin{corollary}
     Let $F$ be the multiperiodic function defined by \eqref{eq-multiperiodic}.
     Suppose that $\beta>1$ is a PV--number and that $f_1, \cdots, f_d$
     are either identically zero or larger than a constant $\delta
     >0$. Suppose further that $M_x^\ell := M(\beta^{\ell-1}
     x)\cdots M(\beta x)M(x)$ has strictly positive entries
     for some integer $\ell>0$.
     Then for any $q \in \mathbb{R}^+$, the following limit
     exists
     \[
             \lim_{T \to \infty} \frac{1}{\log T} \int_{0}^T
             F(x)^q d x.
     \]
\end{corollary}

\begin{proof}
Without loss of generality, we may only consider the subsequence
$T_n= \beta^n$. Since $|G(x)| = \sum_{j=0}^{d-1} |F(x/\beta^j)|$
where $G$ is the solution of the associated vector
equation~\eqref{eq-vector}, we have only to show the existence
of the limit
\[
             \lim_{n \to \infty} \frac{1}{n} \int_{0}^{\beta^n}
             |G(x)|^q d x.
     \]
Making the change of variables $x=\beta^n y$, we are led to prove
the existence of the limit
\[
             \lim_{n \to \infty} \frac{1}{n} \int_{0}^1
             |G(\beta^n x)|^q d x.
     \]
Notice that $G(\beta^{n}x) = M^n_x G(x)$. Notice also that
$G(x)$ has strictly positive entries by the hypothesis on
$M_x^\ell$. So, for the  non negative matrix $M_x^n$ we have
\[
    C^{-1} \| M_x^n\| \le |G(\beta^{n}x)| \le C \|M^n_x\|
         \qquad (\forall x \in [0,1])
\]
for some constant $C>0$. By the proof of the last theorem, $ \log
\int_{0}^1 |G(\beta^n x)|^q d x$ is subadditive.
\end{proof}

\begin{remark} {\rm
   Let $M(x)$ be the matrix defined by~\eqref{Matrix-multip}.
   Let $\tilde{M}$ be the numerical matrix obtained by replacing
   $f_j(x)$ in $M(x)$ by $0$ or $1$ according to $f_j(x) \equiv 0$
   or not. Then $\tilde{M}^\ell>0$ implies $M_x^\ell >0$.
   In particular, $\tilde{M}^d >0$ if $f_j(x)$ are all strictly
   positive.
}
\end{remark}

\begin{example} {\rm Let $\beta >1$ be a PV--number. Let $f_1(x)$ and
$f_2(x)$ be two strictly positive $1$--periodic H\"{o}lder
continuous functions such that $f_1(0) + f_2(0)= 1$. There is a
unique multiperiodic function $F$ defined by
\[
              F(x) = f_1\left(\frac{x}{\beta}\right)F\left(\frac{x}{\beta}\right)
               + f_2\left(\frac{x}{\beta^2}\right)F\left(\frac{x}{\beta^2}\right)
\]
For almost every $x \in \mathbb{R}$, $n^{-1} \log F(\beta^n x)$
has a limit as $n \to \infty$; for any $q \in \mathbb{R}^+$, $(\log
T)^{-1} \int_0^T F(x)^q d x$ has a limit as $T \to \infty$. }
\end{example}

\begin{example} {\rm
Let $\beta>1$ and $a, b \in \mathbb{Z}$. Consider the contractive
transformations on $\mathbb{R}$ defined by
\[
  S_1 x = \frac{x +a}{\beta}, \qquad  S_1 x = \frac{x
  +b}{\beta^2}.
\]
For any $0<p<1$, there exists a unique probability measure $ \mu$
with compact support such that
\begin{equation*}~\label{sim-measure}
         \mu = p \ \mu \circ S_1^{-1} + (1-p)\
          \mu \circ S_2^{-1}.
\end{equation*}
It is a self-similar measure. Its Fourier transform satisfies the
equation
\begin{equation*}
     \widehat{\mu} (x) = f_1(x/\beta) \widehat{\mu} (x/\beta)
         + f_2(x/\beta^2) \widehat{\mu} (x/\beta^2)
\end{equation*}
with $f_1(x) = p e^{2 \pi i a x}$ and $f_2(x) = q e^{2 \pi i b x}$
with $q=1-p$ . This is a special case of the
equation~\eqref{eq-multiperiodic}.
The corresponding matrix defined by~\eqref{Matrix-multip} and its
inverse are respectively equal to
\[
   M(x) =  \left(
    \begin{array} {ll}
             p e^{2 \pi i a x} & q e^{2 \pi i b x/\beta}\\
             1   & 0
    \end{array}
    \right),\quad
     M(x)^{-1} =  q^{-1} e^{- 2\pi i b x/\beta}\left(
    \begin{array} {ll}
               0 &  q e^{2 \pi i b x/\beta} \\
             1   &   -p e^{2 \pi i a x}
    \end{array}
    \right)
\]
 If we take the norm $|v|= \max(|v_1|, |v_2|)$ on $\mathbb{C}^2$,
the operator norms for $M(x)$ and $M(x)^{-1}$ are respectively
$\|M(x)\| = 1$ and $\|M(x)^{-1}\| = \frac{1+ p}{1-p}$. So, when
$\beta $ is a PV--number, under the condition
    $\frac{1+p}{1-p}<\frac{1}{\rho}$, for almost all $x \in \mathbb{R}$
    the following limit exists and does not depend on $x$
    \[
          \lim_{n \to \infty} \frac{1}{n} \log \left( |\widehat{\mu}(\beta^nx)|
          + |\widehat{\mu}(\beta^{n-1}x)| \right).
    \]
    }
\end{example}

\end{document}